\documentclass[a4,12pt]{article}
\usepackage{amssymb}
\catcode`@=11
\def\ps@plain{\let\@mkboth\@gobbletwo
	\def\@oddhead{\rm\hfil\thepage\hfil}
	\def\@oddfoot{}\def\@evenfoot{}\let\@evenhead\@oddhead}
\catcode`@=12 
\normalbaselines
\newcommand{\Z}{\mathbb{Z}}

\newcommand{\R}{\mathbb{R}}
\newcommand{\C}{\mathbb{C}}

\newcommand{\lie}[1]{%
\mathfrak{{#1}}}

\newcommand{\iunit}{%
\sqrt{-1 \,} \, }
\newcommand{\half}{%
\frac{\, 1\,}{\, 2\,}}

\newcommand{\proof}{%
\noindent
{\sc Proof.} \ }
\newcommand{\qed}{%
\hspace*{\fill} Q.E.D.

\bigskip}

\newcommand{\smallspace}{\hspace*{0.8em}}
\newcommand{\medspace}{\hspace*{1.5em}}
\newcommand{\bigspace}{\hspace*{2.5em}}

\newcommand{\0}{\bar{0}}
\newcommand{\1}{\bar{1}}

\newcommand{\sspace}[1]{%
{#1}={#1}_{\0 } \oplus {#1}_{\1 }}
\newcommand{\qroot}{%
\sqrt[4]{-1}}

\newcommand{\ad}{{\rm ad}\;}

\newtheorem{etheorem}{{\sc Theorem}}[section]
\newtheorem{eproposition}[etheorem]{\sc Proposition}

\newtheorem{elemma}[etheorem]{\sc Lemma}
\newtheorem{edefinition}[etheorem]{\sc Definition}

\begin{document}


\begin{center}
{\large \bf Realization of irreducible unitary \\
representations of
{\normalsize $ \lie{osp}(M/N;\R ) $} \large \bf on Fock spaces.
\footnote{
appeared in 
``\textit{Representation theory of Lie groups and Lie algebras}'' (Fuji-Kawaguchiko, 1990), 1–21, World Sci. Publ., River Edge, NJ, 1992.
}
}

\bigskip

Dedicated to Professor Nobuhiko Tatsuuma on his 60th birthday

\bigskip

\bigskip

\small

{\bf By}\\
\bigskip

{\bf Hirotoshi FURUTSU}\\
\medskip
{\it Department of Mathematics, College of Science and Technology,\\
Nihon University\\
Kanda-Surugadai 1-8, Chiyoda, Tokyo 101, JAPAN.}

\bigskip

{\bf Kyo NISHIYAMA}\\
\medskip
{\it Institute of Mathematics, Yoshida College, Kyoto University\\
Sakyo, Kyoto 606, JAPAN.}

\end{center}

\normalsize

\section*{Introduction.}

Recently, Lie superalgebras and their representations have become more important in physics.
They appeared first in the area of elementary particle physics, then they were used in other fields, such as nuclear physics, supergravity, superstring theory, etc. (cf.\  \cite{Kostelecky-Campbell}).
In physics, unitarizability of the representations is very important, since operators in physics have often to be self-adjoint. 
It is notable that many physical applications in supersymmetry deal with Lie superalgebras $ \lie{su}(M/N) $ and $ \lie{osp}(M/N) $. For example in \cite{Bars}, $ \lie{osp}(4/N) $ is used to describe supergravity.

The classification of finite dimensional simple Lie superalgebras was discovered by V.G.~Kac (\cite{Kac1}). 
He also wrote many works on the theory of Lie superalgebras and their finite dimensional representations (\cite{Kac2}, \cite{Kac4}).  
Thereafter, many studies on Lie superalgebras were made in mathematics. 
For example, M.~Scheunert wrote some interesting works already in the 1970's (\cite{Scheunert1}) and F.A.~Berezin also played an important role in this field (\cite{Berezin}).

In the 1970's, most of the studies were of finite dimensional representations (e.g. \cite{Kac4}, \cite{Scheunert.N.R2}).
Two types of irreducible representations of simple Lie superalgebras exist, i.e.\  typical and atypical ones (see \cite{Kac2} for definition). 
For finite dimensional typical representations, similar properties as those of simple Lie algebras hold (\cite{Kac2}, \cite{Kac4}).  
But atypical representations (even finite dimensional ones) are more difficult to study. Many problems are still left to be solved (cf.\  \cite{Gould}).

In recent years, papers dealing with infinite dimensional representations have become more common.  
Now they are studied everywhere, and papers on super unitary (infinite dimensional) representations are also popular (see \cite{Gould-Zhang2}, \cite{Gunaydin}, \cite{Schmitt.et.al}, etc.). It seems interesting that whether a representation is typical or not is not so much correlated with its super unitarity.

Classifications for irreducible super unitary representations were discovered for some low rank basic classical Lie superalgebras. 
For example, for some of the Lie superalgebras of type A or their real forms, classifications were discovered (\cite{Furutsu2}, \cite{Furutsu-Hirai}, \cite{Gould-Zhang1}, \cite{Jakobsen}, etc.).  For an orthosymplectic Lie superalgebra $ \lie{osp}(2/1; \R ) $, a classification was also discovered in \cite{Furutsu-Hirai}.  In \cite{Heidenreich}, irreducible super unitary lowest weight modules of $ \lie{osp}(4/1) $ were classified.  
The complete list of irreducible super unitary representations of $ \lie{osp}(2/2; \R ) $ is obtained in \cite[Th.~4.5]{Nishiyama2}.  For general orthosymplectic algebras $ \lie{osp}(M/N; \R ) \; (M=2m) $, many irreducible super unitary representations were realized explicitly in \cite{Nishiyama2}, using super dual pairs.

For simple Lie algebras, a general classification theorem of irreducible unitary highest weight representations was already given in \cite{EHW} and \cite{Jakobsen3} independently.  
But, for Lie superalgebras, a general classification theorem had been expected for a long time.  Recently, H.P.~Jakobsen made a complete classification of irreducible super unitary highest weight modules (\cite{Jakobsen2}).  
However, his classification only clarifies the parameters of representations and does not provide realizations, character formulas, etc.  In our last paper \cite{Furutsu-Nishiyama}, we classified the irreducible super unitary representations of Lie superalgebras of type $ \lie{su}(p,q/n) $. 
Our method of classification naturally provides realizations of irreducible super unitary representations on Fock spaces.  We think the method is also useful for obtaining character formulas (cf.\  \cite{Nishiyama4}).

In this paper, we classify irreducible super unitary representations of an orthosymplectic Lie superalgebra $ \lie{osp}(M/N; \R ) $ in a similar way as was done in the previous paper.  
Our starting point for classification for $ \lie{su}(p,q/n) $ was based on the notion of super dual pairs in $ \lie{osp}(M/N) $ (cf.\  \cite{Nishiyama2}, \cite{Nishiyama3}).  
Therefore, logically, this paper should have preceded our last paper \cite{Furutsu-Nishiyama}.

In this paper we study the irreducible super unitary representations of a real form of the complex Lie superalgebra $ \lie{osp}(M/N; \C ) \; (N \geq 2) $. 
Note that $ \lie{osp}(M/N; \C ) $ itself has no irreducible super unitary representations except trivial ones.  
The real forms of the Lie superalgebra $ \lie{osp}(M/N; \C ) $ are isomorphic to one of the real Lie superalgebras $ \lie{osp}(M/p,q; \R ) \; (N=p+q, [N/2] \leq p \leq N) $. 
But if $ p $ is not equal to $ N $, $ \lie{osp}(M/p,q; \R ) $ also has no irreducible super unitary representations except trivial ones. So the only real form which has non-trivial super unitary representations is $ \lie{osp}(M/N; \R ) $.

The main result given in this paper is a classification of all the irreducible super unitary representations of $ \lie{osp}(M/N; \R ) \; (N \geq 2) $ that integrate to global representations of $ Sp(M) \times SO(N) $, a Lie group corresponding to the even part (Theorem \ref{thm:cond}).  
As commented above, the classification itself is not new and having been obtained by H.P.~Jakobsen \cite{Jakobsen2} for general highest weight representations in quite different style.  But the construction of integrable unitary representations is new, and we think our method of classification is simpler and gives more information than the one in \cite{Jakobsen2}.

Let us explain the method of classification briefly.  
First of all, we note that if a representation is integrable, then it is admissible (see \S \ref{sec:osp}). Therefore, an irreducible integrable super unitary representation must be a lowest or highest weight representation (see \cite[Prop.~2.3]{Nishiyama2}). 
So what we have to do is to find a necessary and sufficient condition for an irreducible lowest (or highest) weight representation to be super unitary.

{\tolerance=10000
Since an orthosymplectic algebra has a special super unitary representation called oscillator representation (see \cite{Nishiyama1}), we utilize this representation to realize the \linebreak irreducible super unitary representations. 
If we imbed $ \lie{osp}(M/N; \R ) $ into \linebreak
$ \lie{osp}(ML/NL; \R ) $ $ (L \geq 1 ) $, then the oscillator representation of $ \lie{osp}(ML/NL; \R ) $ becomes a super unitary representation of $ \lie{osp}(M/N; \R ) $ through the above imbedding.  
If we decompose it, we can get many irreducible super unitary representations. 
In this paper, the complete decomposition is not carried out, but we construct many primitive vectors for $ \lie{osp}(M/N; \R ) $ explicitly. 
The weights of these vectors are lowest weights for some irreducible super unitary representations. 
This supplies us with a sufficient condition for super unitarity.}

Next we study a necessary condition.  
From the definition of super unitarity, we get an inequality for weights of super unitary representations (Proposition \ref{prop:ncond1}, cf.\  \cite[Prop.~2.2]{Furutsu2} and \cite[Lemma~2.1]{Furutsu-Nishiyama}).  
Note that this inequality also implies that an irreducible admissible super unitary representation must be a lowest (or highest) weight representation.  
It is very important that not only the lowest (or highest) weight but all the weights of the representation must satisfy this condition. 
Take a non-zero lowest weight vector $ v_{\lambda} $ of an irreducible lowest weight super unitary representation, where $ \lambda $ is the lowest weight.  
We choose a special series of weight vectors $ v_1, \cdots , v_{m-1} $ in \S \ref{sec:cond.unitary}.  
If a vector $ v_k $ does not vanish, then its weight must satisfy the above mentioned inequality.  
We translate this condition of the weight of $ v_k $ into that of $ \lambda $, and determine when $ v_k $ does not vanish.  Then we get a necessary condition for super unitarity for the lowest weight $ \lambda $ (Proposition \ref{prop:ncond2}).  

Finally, using this necessary condition for super unitarity, we prove that the above sufficient condition is also necessary, that is, the representations we get by the above method exhaust all the irreducible integrable super unitary representations up to isomorphism.

Our method of classification is quite different from that in \cite{Jakobsen2}. We think our method is simpler in the case of orthosymplectic algebras.  
Moreover it gives us realizations of the representations at the same time. 
On the other hand, we only treat integrable representations, but H.P.~Jakobsen classified all the irreducible super unitary highest weight modules which have continuous classification parameters.

Let us explain each section briefly.  
We introduce some notations for the orthosymplectic algebras $ \lie{osp}(M/N) $ in \S 1.  
We also define the ``admissibility" and ``integrability" for representations of an orthosymplectic Lie superalgebra in \S 1.

In \S 2, we review the oscillator representation of $ \lie{osp}(M/N; \R ) $ which is a super unitary lowest weight representation.  
We imbed $ \lie{osp}(M/N; \R ) $ into $ \lie{osp}(ML/NL; \R ) $ and consider the oscillator representation of $ \lie{osp}(ML/NL; \R ) $ as a super unitary representation of $ \lie{osp}(M/N; \R ) $ through that imbedding.  

In \S 3, we construct a number of $ \Delta^{-} $-primitive vectors in the above representation space (a Fock space), where $ \Delta^{-} $ is the standard negative root system of $ \lie{osp}(M/N; \R ) $ (Proposition \ref{prop:primvec}).  Each of these vectors generates an irreducible super unitary representation of $ \lie{osp}(M/N; \R ) $.  
We calculate the weights of these vectors, then the existence of these weights gives a sufficient condition for super unitarity (Proposition \ref{prop:weight}).

In \S 4, first we introduce a necessary condition for super unitarity which is obtained by the inequality which defines super unitarity (Proposition \ref{prop:ncond1}).  
But this necessary condition is not sufficient. So we must give a stricter necessary condition (Proposition \ref{prop:ncond2}).  
After giving that, we finally show that this condition is also sufficient.  
Thus we get a classification of irreducible integrable super unitary representations of $ \lie{osp}(M/N; \R ) \; (N \geq 2) $.  
This is our main theorem (Theorem \ref{thm:cond}).  
We get realizations of super unitary representations on Fock spaces simultaneously.

\bigskip

The authors whish to express their thanks to Prof.~T.~Hirai for drawing their attention to this field and suggesting that they collaborate.

\section{Orthosymplectic algebras.}
\label{sec:osp}

Let $ \sspace{V} $ be a super space over $ \C $ with dim $ V_{\0} = M $ and dim $ V_{\1} = N $. Throughout in this paper, we put $ M=2m \; (m \in \Z_{\geq 0}) $ and $ N=2n $ or $ 2n+1 \; (n \in \Z_{\geq 0}) $ according as $ N $ is even or odd.
We call a bilinear form $ b $ on $ V $ super skew symmetric if $ b $ satisfies
\[
     b(v,w) = -(-)^{\deg v \deg w} b(w,v),
\]
where $ v, w $ are homogeneous elements in $ V $ and deg $ v $ denotes the degree of $ v $.
We denote an orthosymplectic algebra by
\[
   \lie{osp}(b) = \{ X \in \lie{gl}(V) | b(Xv,w)+(-)^{\deg X \deg v}b(v,Xw)=0                    \mbox{ for } v,w \in V \} ,
\]
where deg $ X $ denotes the degree of $ X \in \lie{gl}(V) $.  
Note that an element $ X $ in the above formula is supposed to be homogeneous.  However, we say that $ X $ satisfies some property (or formula) when homogeneous components of $ X $ satisfy that property (or formula), if there is no confusion.  So $ \lie{osp}(b) $ consists of linear combinations of homogeneous $ X $ which satisfies the above equation.

In this paper, we fix a bilinear form $ b $ which is expressed by a matrix of the form
\[     B= \left[ \begin{array}{cc|c}
           & 1_m &         \\ 
      -1_m &     &         \\ \hline 
           &     & 1_{N}    
      \end{array} \right],
\]
that is $ b(v,w) = {^t}vBw $ for column vectors $ v,w \in V $.
We denote an orthosymplectic algebra $ \lie{osp}(b) $ for the above $ b $ by $ \lie{osp}(M/N; \C ) $ or $ \lie{osp}(M/N) $ simply.
The elements in $ \lie{osp}(M/N; \C ) $ are matrices of the form
\[ \left[ \begin{array}{cc|c}
   A            & B            & P \\        
   C            & - \; ^t \! A & Q \\ \hline 
   - \; ^t \! Q & ^t \! P      & D           
   \end{array} \right] , \]
where $ A \in \lie{gl}(m; \C ) $ , $ B $ and $ C $ are symmetric, $ P $ and $ Q $ are $ m \times N $-matrices, and $ D $ belongs to $ \lie{so}(N) $.  
It is easy to see that real forms of $ \lie{osp}(M/N; \C ) $ are $ \lie{osp}(b_p; \R ) \; ([\frac{N}{2}] \leq p \leq N) $ up to isomorphism, where $ b_p $ is expressed by a matrix
\[     B_p = \left[ \begin{array}{cc|cc}
           & 1_m &       &      \\ 
      -1_m &     &       &      \\ \hline 
           &     & 1_{p} &      \\ 
           &     &       & -1_{N-p} 
      \end{array} \right] .
\]
If $ p \neq N $, these real forms have no admissible unitary representation except for trivial ones (\cite[Prop.~2.3]{Nishiyama2}).
We denote $ \lie{g} = \lie{osp}(b_N; \R ) $ by $ \lie{osp}(M/N; \R ) $.
We will write down the root system of this algebra in order to fix the notations.
This algebra has a compact Cartan subalgebra $ \lie{h} $:
\[
     \lie{h}= \left\{ h= \left[ \begin{array}{c|c} 
             \begin{array}{cc}  0 &  A \\ 
                               -A &  0  \end{array}  &  0 \\ \hline 
                                0 &  B
                         \end{array} \right] \left| 
             \begin{array}{l}
                         A= {\rm diag}(a_1, a_2, \cdots ,a_m), \\
                         B= {\rm diag}(b_1 u, b_2 u, \cdots , b_n u, 1), \\
                         a_i,b_j \in \R
             \end{array} \right. \right\},
\]
where
\[
    u = \left( \begin{array}{cc}
               0 & 1 \\
              -1 & 0
               \end{array} \right) ,
\]
and the figure $ 1 $ in the last place of the sequence in $ B $ appears if and only if $ N $ is odd.
We define $ e_i \in (\lie{h}^{\C})^{\ast} \; (1 \leq i \leq m) $ and $ f_j \in (\lie{h}^{\C})^{\ast} \; (1 \leq j \leq n) $ by putting  
\[     e_i(h)= \iunit  a_i, \hspace{2em}    f_j(h)= \iunit  b_j,  \]
\noindent
for $ h \in \lie{h} $ of the above form. Then, for the case $ N =2n $, roots are given as  
\[ \makebox[13.5cm][l]{$ \Delta _c ^+ =\{ e_i-e_j | 1 \leq i < j \leq m \} \cup 
                   \{ f_i \pm f_j | 1 \leq i < j \leq n \} $} \] 
\hspace*{.52\textwidth}
:the set of positive compact roots, 
\begin{eqnarray*}
& \makebox[6.8cm][l]{$ \Delta _n ^+ =\{ e_i+e_j | 1 \leq i \leq j \leq m \} $ } & 
                 \mbox{:the set of positive non-compact roots,} \\ 
& \makebox[6.8cm][l]{$ \Delta _{\0}^+ = \Delta _c ^+ \cup \Delta _n ^+ $}  & 
                 \mbox{:the set of positive even roots,} \\ 
& \makebox[6.8cm][l]{$ \Delta _{\1}^+ =\{ e_i \pm f_j | 
                               1 \leq i \leq m, 1 \leq j \leq n \} $ } & 
                 \mbox{:the set of positive odd roots,}  \\
& \makebox[6.8cm][l]{$ \Delta^+ = \Delta_{\0}^+ \cup \Delta_{\1}^+ $}  & 
                 \mbox{:the set of positive roots,}  
\end{eqnarray*}
For the case $ N = 2n+1 $, a set of positive roots are similarly given but the set of positive compact roots contains $ \{ f_j | 1 \leq j \leq n \} $
and the set of positive odd roots contains $ \{ e_i | 1 \leq i \leq m \} $:
\[
     \Delta _c ^+ =\{ e_i-e_j | 1 \leq i < j \leq m \} \cup 
                   \{ f_i \pm f_j | 1 \leq i < j \leq n \} \cup
                   \{ f_j | 1 \leq j \leq n \},
\]
\[
     \Delta _{\1}^+ =\{ e_i \pm f_j | 1 \leq i \leq m, 1 \leq j \leq n \} \cup 
                     \{ e_i | 1 \leq i \leq m \}.
\]
We put
\[
     \lie{g}^{\pm }_{\0 } =
     \bigoplus_{\pm \alpha \in \Delta^{+}_{\0 } } \lie{g} _{\alpha },
     \medspace \lie{g}^{\pm }_{\1 } =
     \bigoplus_{\pm \alpha \in \Delta^{+}_{\1 } } \lie{g} _{\alpha },
     \medspace \lie{g}^{\pm} = \lie{g}^{\pm}_{\0} \oplus \lie{g}^{\pm}_{\1},
\]
where $ \lie{g}_{\alpha } $ is a root space in $ \lie{g}^{\C } $ of a root $ \alpha $.

If $ \lambda \in (\lie{h}^{\C})^{\ast} $ is of the form  

\[    \lambda = \sum _{1 \leq i \leq m} \lambda_i e_i + 
                \sum _{1 \leq j \leq n} \mu_j f_j,  \]

\noindent 
then we write $ \lambda =(\lambda _1, \lambda _2 , \cdots , \lambda _m / \; \mu _1 , \mu _2 , \cdots , \mu _n ) $ and call it a {\it coordinate expression} of $ \lambda  $.

\begin{edefinition} \label{def:adm} \rm
A representation $ ( \pi , U ) $ of $ \lie{g} = \lie{osp}(M/N;\R ) $ is called {\it infinitesimally admissible} if the representation $ ( \pi |_{\lie{k}} , U ) $ of $ \lie{k} $ is a direct sum of its irreducible finite dimensional representations with finite multiplicity, where $ \lie{k} $ is a maximal compact subalgebra of $ \lie{g}_{\0 } $.
\end{edefinition}

Note that, to define the notion of admissibility, 
we usually use a compact Lie group $ K $, which has the Lie algebra $ \lie{k} $.   But in the above definition, we use $ \lie{k} $ instead of $ K $. So we add the term ``infinitesimally".  
In the following, we assume that $ \lie{k} $ contains the Cartan subalgebra $ \lie{h} $.

It is known that an irreducible unitary representation of a semisimple Lie group is admissible (see, for example, \cite[Th.~0.3.6]{Vogan1}). Therefore, it is natural to consider infinitesimally admissible representations even for Lie superalgebras.  
Moreover, we only consider ``integrable" representations in this paper.

\begin{edefinition} \label{def:int} \rm 
A representation $ ( \pi , U ) $ of $ \lie{g} = \lie{osp}(M/N;\R ) $ is called {\it integrable} if the representation $ ( \pi |_{\lie{g}_{\0 }} , U ) $ of $ \lie{g}_{\0 } $ is a $ (\lie{g}_{\0}, K) $-module consisting of $ K $-finite vectors of an admissible representation of a Lie group $ G_{\0} = Sp(M) \times SO(N) $, where $ K $ is a maximal compact subgroup of $ G_{\0} $.
\end{edefinition}

It is clear that integrable representations are infinitesimally admissible.  Conversely, let $ ( \pi , U ) $ be an infinitesimally admissible representation of $ \lie{osp}(M/N; \R ) $. Then, it is necessary for $ \pi $ to be integrable that any representation of $ \lie{k} $ in $ \pi $ can be lifted up to a representation of $ K $.  
Therefore any weight $ \lambda $ of $ \pi $ must satisfy the condition
\[
     \lambda_i , \mu_j \in \Z \medspace 
		\mbox{ for } 1 \leq i \leq m, 1 \leq j \leq n,
\]
that is, all $ \lambda _i $'s and $ \mu _j $'s are integers.

We also note that the following proposition holds.

\begin{eproposition}[e.g.\  {\cite[Prop.~2.3]{Nishiyama2}}] \label{prop:admuni}
For $ \lie{osp}(M/N; \R ) $, an irreducible admissible unitary representation is a highest or lowest weight representation.
\end{eproposition}

So we only consider lowest (or highest) weight representations in the rest of this paper.

\section{Oscillator representation.}
\label{sec:osci}

In \cite{Nishiyama1}, we defined a special super unitary representation for $ \lie{osp}(M/N; \R ) $ called the oscillator representation.  Let us review the construction of it briefly for the case $ N =2n+1 $. Note that the case $ N =2n $ is essentially contained in this case (ignore the element $ c $ in the following).

Let $ C^{\C}(r_{\ell} ,c| 1 \leq \ell \leq n) $ be a Clifford algebra over $ \C $ generated by $ \{ r_{\ell} | 1 \leq \ell \leq n \} $ $ \bigcup \{ c \} $ with relations
\[ {r_i}^2 = 1, \hspace{2em}
                 r_i r_j + r_j r_i =0 \; (1 \leq i \neq j \leq n) , \]
\[  c^2 = 1, \hspace{2em}
                 c r_i + r_i c =0 \; (1 \leq i \leq n) . \]
If we set the degree of each generator in the above to $ \1 \in \Z _2 $, then this algebra becomes a superalgebra. That is, let us denote by $ C^{\C}_{\0}(r_{\ell},c| 1 \leq \ell \leq n) $ (resp.\  $ C^{\C}_{\1}(r_{\ell},c| 1 \leq \ell \leq n) $) a subalgebra (resp.\  subspace) of $ C^{\C}(r_{\ell},c| 1 \leq \ell \leq n) $ generated by even (resp.\  odd) products of $ r_j\mbox{'s} $ and $ c $. Then we have
\[    C^{\C}(r_{\ell},c| 1 \leq \ell \leq n) = C^{\C}_{\0}(r_{\ell},c| 1 \leq \ell \leq n) \oplus C^{\C}_{\1}(r_{\ell},c| 1 \leq \ell \leq n),  \]
and this is the decomposition as a $ \Z_2 $-graded algebra.
Now we define a representation space $ \sspace{F} $ of the oscillator representation $ \rho $. Put
\begin{eqnarray*}
\makebox[.6cm][l]{$ F $}   & = & F_{\0} \oplus F_{\1} , \\
\makebox[.6cm][l]{$ F_{\0} $} & = & \C [z_k| 1 \leq k \leq m] \otimes
                           C_{\0}^{\C}(r_{\ell},c | 1 \leq \ell \leq n), \\
\makebox[.6cm][l]{$ F_{\1} $} & = & \C [z_k| 1 \leq k \leq m] \otimes
                              C_{\1}^{\C}(r_{\ell},c | 1 \leq \ell \leq n),
\end{eqnarray*}
where $ \C [z_k| 1 \leq k \leq m] $ means a polynomial algebra generated by $ \{ z_k | 1 \leq k \leq m \} $.

To define the operation $ \rho $ of $ \lie{osp}(M/N; \R ) $ on $ F $, we introduce a representation of a Clifford-Weyl algebra $ C^{\R}(V; b) $ (cf.\  \cite{Tilgner}).
Let $ \sspace{V} $ be a super space of dimension $ (M/N) $ on which $ \lie{osp}(M/N; \R ) $ naturally acts (see \S 1).
And let $ b $ be a super skew symmetric form on $ V $
such that $ \lie{osp}(M/N; \R ) = \lie{osp}(b) $ holds.
Let us choose a basis $ \{ p_k , q_k | 1 \leq k \leq m  \} $ for $ V_{\0} $ such that
\[    b(p_i, q_j)=-b(q_j, p_i)=\delta _{i,j} , \; b(p_i, p_j)=b(q_i, q_j)=0, \]
and an orthogonal basis $ \{ r_{\ell}, s_{\ell}, c | 1 \leq {\ell} \leq n \} $ for $ V_{\1} $ with respect to $ b $ with length $ \sqrt{2} $.  Then there exists a superalgebra $ C^{\R}(V; b) $ over $ \R $ which is generated by $ \{ p_k , q_k | 1 \leq k \leq m  \} \cup \{ r_{\ell} , s_{\ell}, c | 1 \leq \ell \leq n  \} $ with relations
\[  p_i q_j - q_j p_i = \delta_{i,j} , \; r_i s_j + s_j r_i=0 , \;
    r_i r_j + r_j r_i=2\delta_{i,j} , \; s_i s_j + s_j s_i = 2\delta_{i,j} , \]
\[  c^2 = 1 , c r_i + r_i c = 0 , c s_i + s_i c =0 , \]
and all the other pairs of $ p, q, r, s, c $ commute with each other.
$ C^{\R}(V; b) $ can be considered as a Lie superalgebra in the standard way (cf.\  \cite[{\S}1.1]{Kac1}),
and $ \lie{osp}(M/N; \R ) $ can be realized as a sub Lie superalgebra in $ C^{\R}(V; b) $.  
In fact, let $ \lie{L} $ be a subspace generated by the second degree elements of the following form:
\[ \{ xy + (-)^{\deg x \deg y}yx | \; x, y \; 
	\mbox{ are one of the generators }
              p_k , q_k , r_{\ell} , s_{\ell} \mbox{ or } c \} .
\]
Then $ \lie{L} $ becomes a sub Lie superalgebra.  We denote $ m(x,y) = xy + (-)^{\deg x \deg y}yx $.  An operator $ \ad m(x,y) $ preserves $ V \subset C^{\R}(V; b) $ and the bilinear form $ b $, and this adjoint representation gives an isomorphism between $ \lie{L} $ and $ \lie{osp}(M/N; \R ) $.
From now on, we will identify $ \lie{L} $ and $ \lie{osp}(M/N; \R ) $ with each other.  For example, $ \ad m(p_k,r_{\ell}) $
is expressed as a matrix of the form
\[
      \ad m(p_k,r_{\ell}) \longleftrightarrow
      2\sqrt{2} \left[ \begin{array}{cc|c}
        0 &                & E_{k,2\ell -1} \\
          & 0              & 0              \\ \hline
        0 & E_{2\ell -1,k} & 0
      \end{array} \right],
\]
%
%
%
%
%
where $ E_{i,j} $ is an $ m \times N $ or $ N \times m $-matrix with all the elements $ 0 $ except the $ (i,j) $-element, which is $ 1 $.

The oscillator representation $ (\rho , F) $ is actually a representation of the superalgebra $ C^{\R}(V; b) $ given by

\bigskip

\noindent
\renewcommand{\baselinestretch}{2.5} \small \normalsize
\makebox[1.7cm][l]{} $
       \begin{array}{lll}
       \rho (p_{k}) = {\displaystyle \frac{\iunit \qroot}{\sqrt{2}}
        \left( z_{k} - \frac{\partial}{\partial z_k} \right) \otimes 1 }
         & & (1 \leq k \leq m ) , \\
       \rho (q_{k}) =  {\displaystyle \frac{\qroot}{\sqrt{2}} \left( z_k +
         \frac{\partial}{\partial z_k} \right) \otimes 1 }
         & & (1 \leq k \leq m ) , 
       \end{array} $
\renewcommand{\baselinestretch}{1.2} \small \normalsize

\medskip

\noindent
\makebox[2cm][l]{} $
       \rho (r_{\ell})=1 \otimes  r_{\ell} ,  \;
       \rho (s_{\ell})=1 \otimes \iunit r_{\ell} \alpha _{\ell} \smallspace 
	(1 \leq \ell \leq n ), \medspace
       \rho (c)=1 \otimes c, $

\medskip

\noindent
where $ \alpha _{\ell} $ is an automorphism of the Clifford algebra $ C^{\C}(r_{\ell},c| 1 \leq \ell \leq n) $ which sends $ r_k $ to $ (-)^{\delta _{k, \ell}} r_k \; (1 \leq k \leq n) $ and $ c $ to $ c $.
If we restrict $ \rho $ to the sub Lie superalgebra $ \lie{osp}(M/N; \R ) \subset C^{\R}(V; b) $, then $ \rho $ gives a {\it super unitary} representation for $ \lie{osp}(M/N; \R ) $.  For more information on $ \rho $, see \cite{Nishiyama1} and \cite{Nishiyama2}.

Let us consider the following imbedding:
\[
      \iota : \lie{osp}(M/N; \R ) \hookrightarrow \lie{osp}(ML/NL; \R ) 
        \medspace (L \geq 1). 
\]
Here $ \iota $ is given as follows.  Let $ \sspace{V} $ be the super space with the super skew symmetric bilinear form $ b=b_V $ as above and $ W = W_{\0} $ be a usual $ N $-dimensional vector space with a positive definite inner product $ b_W $.  Then a super space $ V \otimes W = V_{\0} \otimes W_{\0} + V_{\1} \otimes W_{\0} $ is endowed with a super skew symmetric bilinear form $ b_{V \otimes W} = b_V \otimes b_W $.  If we consider
\[  \lie{osp}(b_V) = \lie{osp}(M/N; \R ) \]
and
\[  \lie{osp}(b_{V \otimes W}) = \lie{osp}(ML/NL; \R ) , \]
$ \iota $ is given by $ \iota (A) = A \otimes 1_W $ for $ A \in \lie{osp}(M/N; \R ) $.  In the matrix form, this only means 
\[ \iota (A) = \left[ \begin{array}{ccc}
                      a_{11} 1_L & a_{12} 1_L & \cdots \\ 
                      a_{21} 1_L & a_{22} 1_L & \cdots \\ 
                      \vdots     & \vdots     & \ddots    
                      \end{array} \right] \; \mbox{ for }  A = (a_{i,j}) . \]

Let $ (\rho^L , F^L) $ be the oscillator representation of $ \lie{osp}(ML/NL; \R ) $ so that
\[ F^L = \C [z_{i,j} | 1 \leq i \leq m, 1 \leq j \leq L ] \otimes 
	C^{\C}(r_{k,\ell}, c_{\ell} | 1 \leq k \leq n, 1 \leq \ell \leq L ), 
\]
where we use similar notations $ p_{i,k},q_{i,k},r_{\ell,k},s_{\ell,k} $ and $ c_k  \; (1 \leq i \leq m , 1 \leq \ell \leq n , 1 \leq k \leq L) $ as those for $ \lie{osp}(M/N; \R ) $.
We denote $ \tilde{m}(x,y)=xy+(-)^{\deg x \deg y}yx $.  
The successive application of $ \iota $ then $ \rho^L $ gives a {\it super unitary} representation $ \tilde{\rho} = \rho^L \circ \iota $ of $ \lie{osp}(M/N; \R ) $.  In the following subsections, we try to decompose this super unitary representation $ \tilde{\rho} $ of $ \lie{osp}(M/N; \R ) $.  From the definition of super unitarity, there exists a constant $ \varepsilon = \pm 1 $ for each super unitary representation, which is called the associated constant (see \cite{Furutsu-Hirai}). 
Note that, for $ \lie{osp}(M/N; \R ) $, if $ \varepsilon = -1 $ then a super unitary representation must be a lowest weight representation, and if $ \varepsilon = 1 $ then it must be highest (see \cite{Nishiyama2}). 
In this case, since the associated constant of $ \rho $ is $ \varepsilon = -1 $, an irreducible super unitary representation of $ \lie{osp}(M/N; \R ) $ which appears in $ (\tilde{\rho}, F^L) $ is a lowest weight module.  Therefore what we have to do is to find all the $ \Delta^{-} $-primitive vectors for $ \tilde{\rho} $.

Now let us calculate operators for root vectors.  
Let $ X_{\alpha} \; (\alpha \in \Delta ) $ be a non-zero root vector for a root $ \alpha $ of $ \lie{osp}(M/N; \R ) $.  Then up to a non-zero constant multiple, the operators $ \tilde{\rho}(X_{\alpha}) $ are given as follows.

Root vectors for $ \alpha \in \Delta ^- $ in $ \lie{osp}(M/N) $ :

\renewcommand{\baselinestretch}{3} \small \normalsize
\[
\begin{array}{llc}
{\rm (I)} &
\alpha = -(e_s \pm f_{t}) \; (1 \leq s \leq m , 1 \leq t \leq n ); &
{\displaystyle \sum _{k=1} ^L \frac{\partial}{\partial z_{s,k}} 
                  r_{t ,k}(1 \pm \alpha _{t ,k})} ,  \\ 
{\rm (II)} &
\alpha = -(e_s - e_{t}) \; (1\leq s < t \leq m); &
{\displaystyle \sum _{k=1} ^L z_{t, k} \frac{\partial}{\partial z_{s,k}}} , \\ 
{\rm (III)} &
\alpha = -(e_s + e_t) \; (1 \leq s , t \leq m); &
{\displaystyle \sum _{k=1} ^L \frac{\partial}{\partial z_{s,k}} 
                  \frac{\partial}{\partial z_{t,k}}} , \\ 
{\rm (IV)} &
\alpha = -(f_s \pm f_t) \;(1 \leq s < t \leq n );  &
{\displaystyle \sum _{k=1} ^L r_{s ,k}r_{t ,k}
                  (1 + \alpha _{s ,k})(1 \pm \alpha _{t ,k})} , \\
{\rm (V)} &
\alpha = - e_s \; (1 \leq s \leq m); &
{\displaystyle \sum _{k=1} ^L \frac{\partial}{\partial z_{s,k}} c_k } , \\ 
{\rm (VI)} &
\alpha = -f_s \;(1 \leq s \leq n );  &
{\displaystyle \sum _{k=1} ^L r_{s ,k} c_k (1+\alpha _{s ,k})} . 
\end{array} \]
\renewcommand{\baselinestretch}{1.2} \small  \normalsize

All these formulas follow after easy calculations. Here we show only the case (I). 
A root vector $ X_{\alpha} \in \lie{osp}(M/N; \C ) $ of root $ \alpha = -(e_s \pm f_t) $ is expressed by a matrix of the form
\[    \qroot  \left[ \begin{array}{cc|c}
      0 &   & E_{s,2t-1} \mp \iunit E_{s,2t}  \\ 
        & 0 & -\iunit E_{s,2t-1} \mp E_{s,2t} \\ \hline 
      \iunit E_{2t-1,s} \pm E_{2t,s} & E_{2t-1,s} \mp \iunit E_{2t,s} & 0 
      \end{array} \right].
\]
Let us identify $ \lie{L}^{\C} $ with $ \lie{osp}(M/N; \C ) $. Then $ X_{\alpha} $ is identified with an element
\[
    \frac{\qroot}{2\sqrt{2}} \{ m(p_s,r_t) - \iunit m(q_s,r_t)
                    \mp \iunit m(p_s,s_t) \mp m(q_s,s_t) \}.
\]
Therefore we get
\begin{eqnarray*}
\lefteqn{  \rho (X_{\alpha}) = \frac{\qroot}{\sqrt{2}} 
	\{  \rho (p_s) \rho (r_t) - \iunit \rho (q_s) \rho (r_t) 
	\mp  \iunit \rho (p_s) \rho (s_t) 
	\mp  \rho (q_s) \rho (s_t) \} } 	\\
  & = &  \frac{\qroot}{\sqrt{2}} \rho (p_s - \iunit q_s) 
		\rho (r_t \mp \iunit s_t) \\
  & = &  \frac{\qroot}{\sqrt{2}} \left\{
 \frac{\iunit \qroot}{\sqrt{2}} 
	\left( z_s - \frac{\partial}{\partial z_s} \right)
	-\iunit \frac{\qroot}{\sqrt{2}} 
	\left( z_s + \frac{\partial}{\partial z_s} \right) \right\} \otimes \\
  & & \hspace{5cm} \otimes \left\{ \left( r_t \right) 
	\mp \iunit \left( \iunit r_t \alpha _t \right) \right\}  \\
  & = & \frac{\partial}{\partial z_s} r_t (1 \pm \alpha _t). 
\end{eqnarray*}
Next we consider an operator $ \tilde{\rho}(X_{\alpha}) $ for $ \lie{osp}(M/N; \R ) $. Then similarly as above, 
%
%
$ X_{\alpha} $ is identified with an element
\[
    \frac{\qroot}{2\sqrt{2}} 
	\sum_{k=1}^{L} \{ \tilde{m}(p_{s,k},r_{t,k}) 
	- \iunit \tilde{m}(q_{s,k},r_{t,k}) 
	\mp \iunit \tilde{m}(p_{s,k},s_{t,k}) 
	\mp \tilde{m}(q_{s,k},s_{t,k}) \} .
\]
Thus we get
\[
  \tilde{\rho} (X_{\alpha}) = \frac{\qroot}{2\sqrt{2}} 
	\sum_{k=1}^{L} \{ 2 \tilde{\rho} (p_{s,k}) \tilde{\rho} (r_{t,k}) 
	-2 \iunit \tilde{\rho} (q_{s,k}) \tilde{\rho} (r_{t,k})  
	\hspace{3cm}
\]
\[ \hspace{3cm}
	\mp 2 \iunit \tilde{\rho} (p_{s,k}) \tilde{\rho} (s_{t,k}) 
	\mp 2 \tilde{\rho} (q_{s,k}) \tilde{\rho} (s_{t,k}) \} .
\]
Similar calculations lead us to
\[
  \tilde{\rho} (X_{\alpha}) = \sum _{k=1} ^L 
	\frac{\partial}{\partial z_{s,k}} r_{t,k} (1 \pm \alpha _{t,k}).
\]
This completes the discussion of the operators for root vectors.

Finally we show operators for the Cartan subalgebra $ \lie{h} $ of $ \lie{osp}(M/N; \R ) $.  Let $ A_i $ be an element of $ \lie{h} $ which is expressed by a matrix of the form
\[    \left[ \begin{array}{cc|c}
      0 & E_{i,i} &  \\ 
     -E_{i,i} & 0 &  \\ \hline 
        &   & 0 
      \end{array} \right]
      \medspace (1 \leq i \leq m).
\]
Let us identify $ \lie{L} $ with $ \lie{osp}(M/N; \R ) $. Then $ A_i $ is identified with an element
\[
     \frac{1}{4} \{ m(p_i,p_i) + m(q_i,q_i) \} .
\]
Therefore we get
\begin{eqnarray*}
\lefteqn{\rho (A_i) = \frac{1}{4} 
	\{ 2\rho (p_i)\rho (p_i) + 2\rho (q_i)\rho (q_i) \} } 	\\
	& = &  - \frac{\iunit }{4}
		{\left( z_i - \frac{\partial}{\partial z_i} \right)}^2
		+\frac{\iunit }{4}
		{\left( z_i + \frac{\partial}{\partial z_i} \right)}^2
	 =  \iunit \left( z_i \frac{\partial}{\partial z_i}+\half \right) .
\end{eqnarray*}
Then, by the similar arguments as above, we get
a formula of a weight operator for $ \lie{sp} $-part:
\[
    \tilde{\rho} (A_i) = \iunit \sum _{k=1}^L \left( z_{i,k} 
		\frac{\partial}{\partial z_{i,k}} +\half \right) 
		\medspace (1 \leq i \leq m).
\]

Let $ B_j $ be an element of $ \lie{h} $ which is expressed by a matrix of the form:
\[    \left[ \begin{array}{cc|c}
      0 &   &  \\ 
        & 0 &  \\ \hline 
        &   & E_{2j-1,2j} - E_{2j,2j-1} 
      \end{array} \right]
      \medspace (1 \leq j \leq n).
\]
Then $ B_j $ is identified with an element 
$ m(r_j,s_j)/4 \in \lie{L} $.  
Therefore we get
\[
    \rho (B_j) = \frac{1}{2} \rho (r_j) \rho (s_j)
       = \half r_j \iunit r_j \alpha _j = \frac{\iunit}{2} \alpha _j .
\]
Thus we get
a weight operator $ \tilde{\rho}(B_j) $ for $ \lie{so} $-part:
\[
      \tilde{\rho} (B_j) = \frac{\iunit}{2} \sum _{k=1} ^L \alpha_{j,k}                 \medspace (1 \leq j \leq n).
\]

\section{Construction of primitive vectors.}
\label{sec:primvec}

In this section, we will construct primitive vectors for $ \lie{osp} (M/N; \R ) $ in $ F^L $, where $ F^L $ is the representation space of the oscillator representation $ \rho^L $ of $ \lie{osp} (ML/NL; \R ) $. The case of $ N = 2n $ is already treated in \cite{Nishiyama2}. The case of $ N = 2n+1 $ is similar, but we will write down primitive vectors for the sake of completeness.

By definition, $ F^L $ is given by:
\[
     F^L = \C [ z_{s,k} | 1 \leq s \leq m, 1 \leq k \leq L ] \otimes
     C^{ \C } ( r_{t,k}, c_k | 1 \leq t \leq n, 1 \leq k \leq L ).
\]
Put $ L=2l $ or $ 2l+1 $ according as $ L $ is even or odd.  For $ 1 \leq a \leq \min (m,l) $, put
{\renewcommand{\arraystretch}{.5}\small\normalsize
\[
     \Lambda _a = \det (z_{s,2k-1}+ \iunit z_{s,2k} )                                          _{\scriptscriptstyle \begin{array}{l}
\scriptscriptstyle                   m-a <s \leq m \\
\scriptscriptstyle                   1 \leq k \leq a
                   \end{array} },
\]
}
and, for $ 1 \leq b \leq n $ and $ 1 \leq j \leq l $, put
\[
     R_b (j) = \prod_{k=1}^{j} (r_{b,2k-1}+ \iunit r_{b,2k} )                                  \prod_{k'=2j+1}^{L} r_{b,k'},
\]
\[
     C = \prod_{k=1}^{l} (c_{2k-1}+ \iunit c_{2k} ).
\]
For non-negative integers $ i_1, \cdots , i_m $ and $ j_1, \cdots , j_n $, we define
\[
     \Lambda = \Lambda (i_1, \cdots , i_m) = \prod_{a=1}^{m} \Lambda _a ^{i_a},
     \medspace  R = R(j_1, \cdots , j_n) = \prod_{b=1}^{n} R_b (j_b) \cdot C.
\]
We define $ R=C $ if $ N = 1 $.  Note that the order of the product in the expression of $ R $ only causes the multiplication by a sign $ \pm 1 $.  
We get the following result.

\begin{eproposition} \label{prop:primvec}
$ v = \Lambda R $ is a primitive vector for $ \Delta ^{-} $ if 
$ 0 \leq j_1 \leq j_2 \leq \cdots \leq j_n $ and $ i_a = 0 $ for $ a > j_1 $.
\end{eproposition}

\proof
First we will prove that $ \Lambda R $ is killed by negative root vectors which are expressed as linear differential operators, i.e. operators type (I)--(VI) except type (III) (see \S \ref{sec:osci}).  
For these linear differential operators, it is easy to see that we can assume $ \Lambda = \Lambda _a $ ($ a \leq j_1 $).

We easily check that the above vectors are killed by operator (II), since it replaces the $ s $-th row by the $ t $-th row in the determinant $ \Lambda _a $. For the other operators, we note that
\begin{equation} \label{eq:z.prim}
     \frac{\partial}{\partial z_{s,2k}} \Lambda _a                                    = \iunit \frac{\partial}{\partial z_{s,2k-1}} \Lambda _a
\end{equation}
holds for arbitrary $ s $ and $ k $. In fact, both sides represent the same cofactor of $ \Lambda _a $ if $ m-a < s \leq m, 1 \leq k \leq a $.  Otherwise both sides vanish at the same time.
Also
\begin{equation} \label{eq:r.prim}
     r_{b,2k}(1 \pm \alpha _{b,2k}) R_b (j_b)
     = \iunit r_{b,2k-1}(1 \pm \alpha _{b,2k-1}) R_b (j_b)
\end{equation}
for $ k \leq j_b $, and for $ k > 2 j_b $, we have
\begin{equation} \label{eq:r2.prim}
     r_{b,k}(1 + \alpha _{b,k}) R_b (j_b) = 0.
\end{equation}
Similarly $ c_{2k} C = \iunit c_{2k-1} C $ holds for $ k \leq l $.  
Using these formulas, we show that operators (I) and (V) kill the vectors in the lemma.  
For $ k \leq a $, we have
\begin{eqnarray*}
\lefteqn{ \left( \frac{\partial}{\partial z_{s,2k-1}} 
		r_{t,2k-1}(1 \pm \alpha _{t,2k-1})
     + \frac{\partial}{\partial z_{s,2k}} 
		r_{t,2k}(1 \pm \alpha _{s,2k}) \right) \Lambda _a R }    \\
  & &   = \left( \frac{\partial}{\partial z_{s,2k-1}} \Lambda _a \right)
       r_{t,2k-1}(1 \pm \alpha _{t,2k-1}) R
     + \iunit \left( \frac{\partial}{\partial z_{s,2k-1}} \Lambda _a \right)
       r_{t,2k}(1 \pm \alpha _{t,2k}) R
\end{eqnarray*}
from equation (\ref{eq:z.prim}). Since $ k \leq a \leq j_1 \leq j_t $, we can apply (\ref{eq:r.prim}) to the above formula, and one can conclude it vanishes.  
For $ k > a $, we get $ \partial \Lambda _a / \partial z = 0 $. Therefore operator (I) kills the vectors.  
For operator (V), we can proceed similarly.

Let us consider operators (IV) and (VI).
For $ k \leq j_s $, we have
\[
     ( r_{s,2k-1}(1 + \alpha _{s,2k-1}) r_{t,2k-1}(1 - \alpha _{t,2k-1})
     + r_{s,2k}(1 + \alpha _{s,2k}) r_{t,2k}(1 - \alpha _{t,2k}) )
     R_s (j_s) R_t (j_t)
\]
\[
     = (-)^{L-j_s} r_{s,2k-1}(1 + \alpha _{s,2k-1}) R_s (j_s)
       r_{t,2k-1}(1 - \alpha _{t,2k-1}) R_t (j_k)
	\hspace{2cm}
\]
\[
	\hspace{2cm}
     - (-)^{L-j_s} r_{s,2k-1}(1 + \alpha _{s,2k-1}) R_s (j_s)
       r_{t,2k-1}(1 - \alpha _{t,2k-1}) R_t (j_t) =0,
\]
from equation (\ref{eq:r.prim}). In the above formula, we have calculated only the essential factor of $ \Lambda _a R $ for operator (IV).
For $ k > 2j_s $, we get
$ r_{s,k}(1 + \alpha _{s,k}) R_s (j_s) =0 $ from equation (\ref{eq:r2.prim}). Thus operator (IV) kills $ \Lambda _a R $.  
For operator (VI), we can argue similarly.

Now we consider operator (III), which is a second degree differential operator.
For $ k \leq 2l $, we have
\begin{eqnarray*}
    \lefteqn{ \left( \frac{\partial}{\partial z_{s,2k-1}} 
		\frac{\partial}{\partial z_{t,2k-1}}
     + \frac{\partial}{\partial z_{s,2k}} 
		\frac{\partial}{\partial z_{t,2k}}\right)
     \Lambda _a \Lambda _b  \bigspace }      \\
   & &     = \frac{\partial}{\partial z_{s,2k-1}} \Lambda _a
             \frac{\partial}{\partial z_{t,2k-1}} \Lambda _b
           + \frac{\partial}{\partial z_{t,2k-1}} \Lambda _a
             \frac{\partial}{\partial z_{s,2k-1}} \Lambda _b  \\
   & &     \hspace{2cm}
	   - \frac{\partial}{\partial z_{s,2k-1}} \Lambda _a
             \frac{\partial}{\partial z_{t,2k-1}} \Lambda _b
           - \frac{\partial}{\partial z_{t,2k-1}} \Lambda _a
             \frac{\partial}{\partial z_{s,2k-1}} \Lambda _b =0,
\end{eqnarray*}
from equation (\ref{eq:z.prim}). For $ k = 2l+1 $, the operator certainly gives $ 0 $. So operator (III) kills $ \Lambda _a \Lambda _b $. For a general $ \Lambda $, the proof is similar.
\qed

Let us give the weight of this primitive vector $ v $.

\begin{eproposition} \label{prop:weight}
The weight of $ v = \Lambda (i_1, \cdots , i_m) R(j_1, \cdots , j_n) $ is given by $ \lambda = ( \lambda _1, \cdots , \lambda _m / \mu _1, \cdots , \mu _n ) $ with
\[
     \lambda _s = \sum_{a=m-s+1}^{m} i_a + \frac{L}{2} \medspace 
     ( 1 \leq s \leq m ), \medspace
%
     \mu _t = j_t - \frac{L}{2} \medspace ( 1 \leq t \leq n ).
\]
\end{eproposition}

\proof
Note that weight operators for $ \lie{sp} $-part are
\[
      \iunit \sum _{k=1}^L \left( z_{s,k} 
		\frac{\partial}{\partial z_{s,k}} +\half \right) 
		\medspace (1 \leq s \leq m).
\]
Because the degree of $ \Lambda _a $ with respect to $ { z_{s,1}, \cdots , z_{s,L} } $ is $ 0 $ for $ 1 \leq s \leq m-a $ and $ 1 $ for $ m-a < s \leq m $, we obtain
\[
     \sum _{k=1}^L z_{s,k} \frac{\partial}{\partial z_{s,k}} \Lambda _a
     = \left\{ \begin{array}{ll}
               0 & \mbox{for $ 1 \leq s \leq m-a $,} \\
               1 & \mbox{for $ m-a < s \leq m $.}
               \end{array} \right.
\]
Therefore we have
\[
     \lambda _s = \sum_{a=m-s+1}^{m} i_a + \frac{L}{2} \medspace 
     ( 1 \leq s \leq m ).
\]

Since weight operators for $ \lie{so} $-part are
\[
\frac{\iunit}{2} \sum _{k=1} ^L \alpha_{t,k} \medspace (1 \leq t \leq n),
\]
we get
\[
     \mu _t = j_t - \frac{L}{2} \medspace ( 1 \leq t \leq n ),
\]
because $ (\alpha _{t,2k-1} + \alpha _{t,2k})(r_{t,2k-1} + \iunit r_{t,2k}) = 0 $.
\qed

\section{Conditions for super unitarity.}
\label{sec:cond.unitary}

In this section, we study necessary and sufficient conditions for super unitarizability of irreducible lowest weight representations of $ \lie{osp}(M/N;\R ) $ for $ N \geq 2 $. 

\begin{eproposition} \label{prop:ncond1}
Let $ ( \pi , V ) $ be an irreducible super unitary lowest weight representation of $ \lie{osp}(M/N;\R ) \; (M = 2m $ and $ N = 2n $ or $ 2n+1 ) $. Then any weight $ \lambda $ of $ ( \pi , V ) $ must satisfy
\[
   |\mu _t| \leq \lambda _s  \medspace ( 1 \leq s \leq m , 1 \leq t \leq n ).
\]
\end{eproposition}

For the proof, see the proof of Proposition 2.3 in \cite{Nishiyama2}. 

From this proposition and the property of $ \lie{g}_{\0} $-lowest weights, the lowest weight $ \lambda $ of an irreducible super unitary lowest weight representation of $ \lie{osp}(M/N;\R ) $ must satisfy the following condition (I):
\[
{\rm (I)} \medspace \left\{ \begin{array}{l}
     -\mu _n \leq \cdots \leq -\mu _1 \leq \lambda _1                                \leq \cdots \leq \lambda _m,\\
     | \mu _n | \leq -\mu _{n-1}  \mbox{ for } N = 2n,                                \smallspace \mu _n \leq 0  \mbox{ for } N = 2n+1.
     \end{array} \right.
\]
Note that $ \mu _1, \cdots , \mu _{n-1} $ are all non-positive integers or all non-positive half integers.

Now we restrict ourselves to the case that $ N $ is strictly greater than $ 1 $. Let $ ( \pi , V ) $ be an irreducible lowest weight representation of $ \lie{osp}(M/N;\R ) \;  (  N \geq 2  ) $ with the lowest weight $ \lambda \in (\lie{h}^{\C})^{\ast} $ and let $ v_0 $ be a non-zero lowest weight vector. We put
\begin{equation} \label{eq:vkdef}
v_k = X_{m-k+1} X_{m-k+2} \cdots X_m v_0 \smallspace \mbox{ for } k = 1,2, \cdots ,m-1,
\end{equation}
where $ X_j $ is a non-zero root vector for $ \beta _j = e_j - f_1 $.

We find that
$ X_{\alpha } v_k = 0 $ for a negative even root vector $ X_{\alpha} \in \lie{g} _{\alpha} $ ($ \alpha \in \Delta _{ \0 }^{-} $).
In fact, we have
\[
     X_{\alpha} v_k = \sum_{i=1}^{k} X_{m-k+1} \cdots X_{m-k+i-1} [ X_{\alpha} ,     X_{m-k+i} ] X_{m-k+i+1} \cdots X_m v_0
\]
\[
     \hspace{5cm} +X_{m-k+1} \cdots X_m X_{\alpha} v_0 ,
\]
where the second term of the right hand side vanishes since $ v_0 $ is lowest. Each bracket 
$ [ X_{\alpha} , X_j ] $ of the first term is in $ \lie{g} _{ \beta _l } ( j<l ) $ or in $ \lie{g} _{-e_s - f_1} $. 
Therefore each member of the first term vanishes since $ (X_l)^{2} = 0 $ or
$ [ X_{-e_s - f_1} , X_t ] = 0 $ respectively. Thus the vector $ v_k $ is primitive for the even part.

\begin{elemma} \label{lemma:3eq}

The following three conditions are mutually equivalent.

\noindent
$ \begin{array}{lc}
     {\rm (a)} &  v_k \neq 0 , \\
     {\rm (b)} &  X_{- \beta _m } X_{- \beta _{m-1} }                                             \cdots X_{- \beta _{m-k+1} } v_k \neq 0 , \\
     {\rm (c)} &  \prod_{l=1}^{k} ( \lambda _{m-l+k} + \mu _1 -l+1 ) \neq 0. 
\end{array} $
\end{elemma}

\proof
First we show conditions (b) and (c) are equivalent. We have
\[
X_{- \beta _{m-k+1} } v_k = X_{- \beta _{m-k+1} } X_{ \beta _{m-k+1} } v_{k-1}
	\hspace{5cm}
\]
\[
 = [ X_{ \beta _{m-k+1} } , X_{- \beta _{m-k+1} } ] v_{k-1}
 + \sum_{l=1}^{k-1} (-)^{k-l} X_{ \beta _{m-k+1} } \cdots X_{ \beta _{m-l} }
 [ X_{ \beta _{m-l+1} } , X_{- \beta _{m-k+1} } ] v_{l-1} 
\]
\[
 + (-)^{k} X_{ \beta _{m-k+1} } \cdots X_{ \beta _{m} } X_{- \beta _{m-k+1} } v_0
\]
\[
 = H_{ \beta _{m-k+1} } v_{k-1}
 + \sum_{l=1}^{k-1} (-)^{k-l} X_{ \beta _{m-k+1} } \cdots X_{ \beta _{m-l} }
 X_{e_{m-l+1} - e_{m-k+1} } v_{l-1},
\]
where $ X_{e_{m-l+1} - e_{m-k+1} } v_{l-1} $ vanishes since $ e_{m-l+1} - e_{m-k+1} $ is in $ \Delta _{\0}^{-} $ and $ v_{l-1} $ is primitive for the even part. So the above formula becomes
\[
    ( \lambda + e_{m-k+2} + \cdots + e_m -(k-1) f_1 )                               ( H_{ \beta _{m-k+1} }) v_{k-1}
     = ( \lambda _{m-k+1} + \mu _1 -k+1) v_{k-1}.
\]
Using induction on $ k $, we get
\[
 X_{- \beta _m } X_{- \beta _{m-1} } \cdots X_{- \beta _{m-k+1} } v_k
 = \prod_{l=1}^{k}  ( \lambda _{m-l+k} + \mu _1 -l+1 ) v_0 .
\]
Thus (b) and (c) are equivalent.

It is clear that (b) $ \Rightarrow $ (a). We will show (a) $ \Rightarrow $ (b).
Suppose $ v_k \neq 0 $. Then
\begin{equation} \label{eq:vkuniv}
      v_0 \in U( \lie{g} ) v_k ,
\end{equation}
since $ V $ is irreducible. Therefore we can write
\[
      v_0 = \sum Y^{+} Y_{\1}^{-} Y_{\0}^{-} v_k,
\]
where $ Y^{+}, Y_{\1}^{-} $ and $ Y_{\0}^{-} $ are monomials of root vectors in $ \lie{g} ^{+}, \lie{g} _{ \1 }^{-} $ and $ \lie{g} _{ \0 }^{-} $ respectively.
If $ Y^{+} $ is not a scalar, the weight of $ Y_{ \1 }^{-} Y_{ \0 }^{-} v_k $ is lower than $ \lambda $.
Thus we may assume $ Y^{+}=1 $. On the other hand, since $ v_k $ is primitive for the even part, we can also assume $ Y_{ \0 }^{-}=1 $ and we get
\[
      v_0 = \sum Y_{ \1 }^{-} v_k.
\]
The weight of $ Y_{ \1 }^{-} $ is the difference of the weight of $ v_0 $ and that of $ v_k $, i.e. $ k f_1 -(e_{m-k+1} + \cdots + e_m) $. Therefore $ Y_{ \1 }^{-} $ must be of the following form:
\[
    Y_{ \1 }^{-} = c X_{- \beta _m } X_{- \beta _{m-1} }                            \cdots X_{- \beta _{m-k+1} } \hspace{1cm} ( c \in \C ),
\]
and we get
\[
    v_0 = c X_{- \beta _m } X_{- \beta _{m-1}} \cdots X_{- \beta _{m-k+1}} v_k.
\]
This shows that
$ X_{- \beta _m } X_{- \beta _{m-1}} \cdots X_{- \beta _{m-k+1}} v_k \neq 0 $.
\qed

\begin{eproposition} \label{prop:ncond2} 
Let $ ( \pi , V ) $ be an irreducible lowest weight representation of $ \lie{osp}(M/N;\R) \;  (  N \geq 2  ) $ with the lowest weight $ \lambda \in (\lie{h}^{\C})^{\ast} $. Assume that $ (\pi , V) $ is admissible.  Then the following condition is necessary for $ ( \pi , V ) $ to be super unitary.

\noindent
The lowest weight $ \lambda $ satisfies
\begin{equation} \label{eq:2'}
    \lambda _1 + \mu _1 \in \{ d,d+1, \cdots , m-1 \} \cup [m-1, \infty),
\end{equation}
where  $ d = \# \{ 1 \leq k \leq m | \lambda _k > \lambda _1 \}. $
\end{eproposition}

\proof
Let $ ( \pi , V ) $ be super unitary. From condition (I), we have $ \lambda _1 + \mu _1 \geq 0 $.  
If $ m = 1 $, then condition (\ref{eq:2'}) only means $ \lambda _1 + \mu _1 \geq 0 $, and there is nothing to prove.  
We consider the case $ m \geq 2 $.  
If $ \lambda _1 + \mu _1 \geq m-1 $ then the condition trivially holds. So we suppose $ k \leq \lambda _1 + \mu _1 < k+1 $ for $ k = 0,1, \cdots , m-2 $.
Note that $ v_{k+1} $ vanishes. In fact, if $ v_{k+1} \neq 0 $, its weight is $ \lambda + ( e_{m-k+1} + \cdots + e_m ) - k f_1 $. Applying Proposition \ref{prop:ncond1}, we get $ \lambda _1 + \mu _1 \geq k+1 $, a contradiction.
So, from Lemma \ref{lemma:3eq}, there exists $ 1 \leq l \leq k+1 $ such that
\[
   \lambda _{m-l+1} + \mu _{1} = l-1.
\]
For $ l \leq k $, the left hand side of the above equation is greater than $ \lambda _1 + \mu _1 \geq k $ and the right hand side is less than $ k-1 $. Thus the above equation must hold for $ l = k+1 $ and $ \lambda _{m-k} + \mu _1 = k $.
Since $ k- \mu _1 \leq \lambda _1 \leq \cdots \leq \lambda _{m-k} = k- \mu _1 $, the equations must hold: $ \lambda _1 = \cdots = \lambda _{m-k} = k- \mu _1 $,
hence $ d \leq k $. So we get $ \lambda _1 + \mu _1 = k \geq d $.
\qed

From these propositions and Proposition \ref{prop:weight}, we get the following result.

\begin{etheorem} \label{thm:cond} 
Let $ (\pi , V) $ be an integrable irreducible super unitary representation of $ \lie{osp}(M/N;\R) \; (N \geq 2) $, which is necessarily a lowest or highest weight representation.  

{\rm (i)}  If $ (\pi , V) $ is a lowest weight representation, then its lowest weight $ \lambda= $ \linebreak
$ (\lambda _1, \cdots , \lambda _m / $ $ \mu _1, \cdots , \mu _n) $ $ \in (\lie{h}^{\C})^{\ast} $ must satisfy conditions {\rm (I)} and {\rm (II):}

\medskip

\noindent
$ \displaystyle 
{\rm (I)} \medspace \left\{ \begin{array}{l}
     -\mu _n \leq \cdots \leq -\mu _1 \leq \lambda _1                                \leq \cdots \leq \lambda _m,\\
     	|\mu_n| \leq -\mu _{n-1} \smallspace \mbox{ for } N = 2n, \medspace 
	\mu _n \leq 0  \smallspace \mbox{ for } N = 2n+1,\\
     	\lambda _i , \mu _j \in \Z  \smallspace \mbox{ for all } i, j,
\end{array} \right.
$

\smallskip

\noindent
$ \displaystyle 
{\rm (II)} \medspace
    \lambda _1 + \mu _1 \geq d
\medspace \mbox{ where } d = \# \{ 1 \leq k \leq m| \lambda _k > \lambda _1 \}.
$

\medskip

Conversely, let $ (\pi , V) $ be an irreducible lowest weight representation of \linebreak
$ \lie{osp}(M/N;\R ) $ $ (N \geq 2) $ with the lowest weight $ \lambda = ( \lambda _1, \cdots , \lambda _m / $ $ \mu _1, \cdots , \mu _n ) \in (\lie{h}^{\C})^{\ast} $ which satisfies the above conditions {\rm (I)} and {\rm (II)}, then $ ( \pi , V ) $ is super unitary.

{\rm (ii)} If $ (\pi, V) $ is a highest weight representation, then its highest weight $ \lambda= $ \linebreak
$ (\lambda _1, \cdots , \lambda _m / $ $ \mu _1, \cdots , \mu _n) $ $ \in (\lie{h}^{\C})^{\ast} $ must satisfy conditions ${\rm (I')}$ and ${\rm (II')}${\rm :}

\medskip

\noindent
$ \displaystyle 
{\rm (I')} \medspace \left\{ \begin{array}{l}
     -\mu _n \geq \cdots \geq -\mu _1 \geq \lambda _1                                \geq \cdots \geq \lambda _m,\\
     	| \mu _n | \leq \mu _{n-1}  \smallspace \mbox{ for } N = 2n, \medspace 
	\mu _n \geq 0  \smallspace \mbox{ for } N = 2n+1,\\
     	\lambda _i , \mu _j \in \Z  \smallspace \mbox{ for all } i, j,
\end{array} \right.
$

\smallskip

\noindent
$ \displaystyle 
{\rm (II')} \medspace
    \lambda _1 + \mu _1 \leq -d
\medspace \mbox{ where } d = \# \{ 1 \leq k \leq m| \lambda _k < \lambda _1 \}.
$

\medskip

Conversely, let $ (\pi , V) $ be an irreducible highest weight representation of \linebreak
$ \lie{osp}(M/N;\R ) $ $ (N \geq 2) $ with the highest weight $ \lambda = ( \lambda _1, \cdots , \lambda _m / $ $ \mu _1, \cdots , \mu _n ) \in (\lie{h}^{\C})^{\ast} $ which satisfies the above conditions ${\rm (I')}$ and ${\rm (II')}$, then $ ( \pi , V ) $ is super unitary.
\end{etheorem}

\proof
If $ \pi $ is integrable, then it is admissible. Thus it must be a lowest or highest weight representation (Proposition \ref{prop:admuni}).  We will prove the statements only for lowest weight representations.  The proof for highest weight representations is similar.  
Note that integrability forces $ \lambda_i $'s and $ \mu_j $'s to be integers, as we mentioned in \S \ref{sec:osp}.  
Let $ \pi $ be an integral lowest weight representation. Then conditions (I) and (II) follow from Propositions \ref{prop:ncond1} and \ref{prop:ncond2}.

If $ \lambda _1 = 0 $, then all $ \mu _j $ vanish because of condition (I) and, from condition (II), we have $ d = 0 $.
Thus all $ \lambda _i $ also vanish. Therefore $ \lambda $ is the weight of the trivial representation which is super unitary.

If $ \lambda _1 \neq 0 $, then let us put $ L = 2 \lambda _1 $ , $ i_k = \lambda _{m-k+1} - \lambda _{m-k} $ for $ 1 \leq k \leq m-1 $ , $ i_m = 0 $ and $ j_b = \lambda _1 + \mu _b \smallspace ( 1 \leq b \leq n ) $.
Then all $ i_k $'s and $ j_b $'s become integers and,
from conditions (I) and (II), these $ L,{i_a},{j_b} $ satisfy the conditions in Proposition \ref{prop:primvec}. From Proposition \ref{prop:weight}, the weight of the corresponding primitive vector $ v $ is $ \lambda $ itself.
Therefore $ \lambda $ is the lowest weight of a super unitary representation.
\hfill Q.E.D.


\begin{thebibliography}{10}

\small

\bibitem{Bars}
I. Bars.
\newblock Supergroups and superalgebras in physics.
\newblock {\it Physica}, {\bf 15D}(1985), 42--64.

\bibitem{Berezin}
F.A. Berezin.
\newblock {\it Introduction to Superanalysis}.
\newblock Reidel, 1987.

\bibitem{EHW}
T. Enright, R. Howe, and N. Wallach.
\newblock A classification of unitary highest weight modules.
\newblock In {\it Representation Theory of Reductive Groups}, pages~97--143,
  Birkh\"{a}user, 1983.

\bibitem{Furutsu2}
H. Furutsu.
\newblock Representations of {L}ie superalgebras, {II}. {U}nitary
  representations of {L}ie superalgebras of type $ {A}(n, 0) $.
\newblock {\it J. Math. Kyoto Univ.}, {\bf 29}(1989), 671--687.

\bibitem{Furutsu-Hirai}
H. Furutsu and T. Hirai.
\newblock Representations of {L}ie superalgebras, {I}. {E}xtensions of
  representations of the even part.
\newblock {\it J. Math. Kyoto Univ.}, {\bf 28}(1988), 695--749.

\bibitem{Furutsu-Nishiyama}
H. Furutsu and K. Nishiyama.
\newblock {Classification of irreducible super unitary representations of $
  \lie{su}(p,q/n) $}.
\newblock to appear in {\it Comm. in Math. Phys.}

\bibitem{Gould}
M.D. Gould.
\newblock Atypical representations for {type-I Lie} superalgebras.
\newblock {\it J. Phys. A., Math. Gen.}, {\bf 22}(1989), 1209--1221.

\bibitem{Gould-Zhang1}
M.D. Gould and R.B. Zhang.
\newblock Classification of all star and grade star irreps of $ gl(n|1) $.
\newblock {\it J. Math. Phys.}, {\bf 31}(1990), 1524--1534.

\bibitem{Gould-Zhang2}
M.D. Gould and R.B. Zhang.
\newblock Classification of all star and grade star irreps of $ gl(m|n) $.
\newblock {\it J. Math. Phys.}, {\bf 31}(1990), 2552--2559.

\bibitem{Gunaydin}
{M. G\"{u}naydin.}
\newblock {Unitary highest weight representations of non-compact supergroups}.
\newblock {\it J. Math. Phys.}, {\bf 29}(1988), 1275--1282.

\bibitem{Heidenreich}
W. Heidenreich.
\newblock {All linear unitary irreducible representations of de Sitter
  supersymmetry with positive energy}.
\newblock {\it Physics Letters}, {\bf 110B}(1982), 461--464.

\bibitem{Jakobsen3}
H.P. Jakobsen.
\newblock Hermitian symmetric spaces and their unitary highest weight modules.
\newblock {\it J. Funct. Anal.}, {\bf 52}(1983), 385--412.

\bibitem{Jakobsen}
H.P. Jakobsen.
\newblock {On the range of unitarity for highest weight representations of
  classical Lie superalgebras}.
\newblock In {\it Topological and geometrical method in field theory, edited by
  J.Hietarinta and J.Westerholm}, pages~103--109, World Scientific, 1986.

\bibitem{Jakobsen2}
H.P. Jakobsen.
\newblock {The full set of unitarizable highest weight modules of basic
  classical Lie superalgebras}.
\newblock K{\o}benhagens Universitet Matematisk Institut, Preprint 18, August
  1989.

\bibitem{Kac1}
V.G. Kac.
\newblock Lie superalgebras.
\newblock {\it Adv. in Math.}, {\bf 26}(1977), 8--96.

\bibitem{Kac2}
V.G. Kac.
\newblock Representations of classical {L}ie superalgebras.
\newblock {\it LNM}, {\bf 676}(1977), 597--626.

\bibitem{Kac4}
V.G. Kac.
\newblock Characters of typical representations of classical {Lie}
  superalgebras.
\newblock {\it Comm. in Alg.}, {\bf 5}(1977), 889--897.

\bibitem{Kostelecky-Campbell}
V.A. Kosteleck\'y and D.K. Campbell.
\newblock Introduction and overview.
\newblock {\it Physica}, {\bf 15D}(1985), 3--21.

\bibitem{Nishiyama1}
K. Nishiyama.
\newblock Oscillator representations for orthosymplectic algebras.
\newblock {\it J. Alg.}, {\bf 129} (1990), 231--262.

\bibitem{Nishiyama2}
K. Nishiyama.
\newblock Super dual pairs and unitary highest weight modules of
  orthosymplectic algebras.
\newblock to appear in {\it Adv. Math.}

\bibitem{Nishiyama3}
K. Nishiyama.
\newblock {Decomposing oscillator representations of $ \lie{osp}(2n/n; \R ) $
  by a super dual pair $ \lie{osp}(2/1;\R ) $ $ \times $$ \lie{so}(n) $}.
\newblock to appear in {\it Comp. Math.}

\bibitem{Nishiyama4}
K. Nishiyama.
\newblock Characters and super-characters of discrete series representations
  for orthosymplectic {L}ie superalgebras.
\newblock {\it J. Alg.}, {\bf 141} (1991), 399--419.

\bibitem{Scheunert1}
M. Scheunert.
\newblock {\it {The theory of Lie superalgebras}}.
\newblock {\it LNM No.\ 716}, Springer, 1979.

\bibitem{Scheunert.N.R2}
M. Scheunert, W. Nahm, and V. Rittenberg.
\newblock Irreducible representations of the $ osp(2,1) $ and $ spl(2,1) $
  graded {Lie} algebras.
\newblock {\it J. Math. Phys.}, {\bf 18}(1977), 155--162.

\bibitem{Schmitt.et.al}
H.A. Schmitt and B.R. Halse, P.~Barrett.
\newblock Positive discrete series representations of the noncompact
  superalgebra $ osp(4/2,\R ) $.
\newblock {\it J. Math. Phys.}, {\bf 30}(1989), 2714--2720.

\bibitem{Tilgner}
H. Tilgner.
\newblock Graded generalizations of {Weyl-} and {Clifford} algebras.
\newblock {\it J. Pure and Appl. Alg.}, {\bf 10}(1977), 163--168.

\bibitem{Vogan1}
D.A. Vogan.
\newblock {\it Representations of real reductive Lie groups}.
\newblock Birkh{\"{a}}user, 1981.

\end{thebibliography}

\end{document}